\documentclass[10pt]{article}

\usepackage[a4paper,margin=2.9cm]{geometry}
\usepackage[english]{babel}
\usepackage{parskip}
\usepackage{amssymb}
\usepackage{amsmath, amsthm, amssymb}
\usepackage{enumerate}
\usepackage{graphicx}
\usepackage{bbm} 
\usepackage{mathrsfs} 
\usepackage{subcaption}
\usepackage{hyperref}

\makeatletter
\renewcommand\section{\@startsection{section}{1}{\z@}%
                                  {-3.5ex \@plus -1ex \@minus -.2ex}%
                                  {2.3ex \@plus.2ex}%
                                  {\normalfont\Large\bfseries}}
\makeatother

\newtheoremstyle{examplestyle}
  {4mm}
  {4mm}
  {\slshape}
  {0pt}
  {\bfseries}
  {\newline}
  {0mm}
  {}

\theoremstyle{examplestyle}
  \newtheorem{Theorem}{Theorem}              

  \newtheorem{Lemma}[Theorem]{Lemma}
  \newtheorem*{Def*}{Definition of dimensions}
    \newtheorem*{Deff*}{Definition}
 \newtheorem*{Defs*}{Auxiliary definitions}
  
  \newtheorem*{Remark*}{Remark}

\begin{document}

\begin{center}
\textbf{\Large The dimension of the Incipient Infinite Cluster.}
\end{center}

\bigskip

\begin{center}
 \textbf{W.P.S. Cames van Batenburg} \\
Radboud University Nijmegen, The Netherlands.\\ E-mail: w.camesvanbatenburg@math.ru.nl\\
\end{center}
\begin{center}
  (\today) \\
\end{center}

\bigskip 


\begin{abstract}
\noindent
We study the Incipient Infinite Cluster (IIC) of high-dimensional bond percolation on $\mathbb{Z}^d$. We prove that the mass dimension of IIC almost surely equals $4$ and the volume growth exponent of IIC almost surely equals $2$.

\end{abstract}


\section{Introduction}

Consider critical nearest-neighbour percolation on $\mathbb{Z}^d$. The \textit{Incipient Infinite Cluster} (IIC) is a random infinite subset of $\mathbb{Z}^d$ which intuitively can be viewed as the critical cluster of the origin, conditioned to be infinitely large. This conditioning induces a new probability measure $\mathbb{P}_{\text{IIC}}$. We study the IIC in \textit{high dimensions} $d$ (see below for formal definitions) and in particular we identify the typical size of IIC under $\mathbb{P}_{\text{IIC}}$. In order to sensibly determine the size of the IIC we use the concepts of \textit{mass dimension} $d_m(A)$ of a subset $A \subset \mathbb{Z}^d$ and the \textit{volume growth exponent} $d_f(G)$ of an infinite connected graph $G$. The former measures the IIC with respect to the (extrinsic) distance of the space $\mathbb{Z}^d$ in which IIC is embedded, while the latter measures the induced graph of IIC with respect to (intrinsic) graph distance. We prove that the mass dimension of IIC is 4 and the volume growth exponent of the graph of IIC is 2, $\mathbb{P}_{\text{IIC}}$-almost surely. See Theorems $1$ and $2$ below. Theorem $1$ gives an explicit and rigorous foundation for the intuition that for high $d$ the IIC is a $4$-dimensional object, a conjecture of physicists going back at least $30$ years \cite{Physics}\cite{Physics2}. 



\subsection{Critical high-dimensional bond percolation}




Let $G=(\mathbb{Z}^d,E)$ be a graph and fix a parameter $p \in [0,1]$. We focus on the case of \textit{nearest-neighbour} bond percolation, meaning that $(x,y) \in E \Leftrightarrow \|x-y\|_{1} = 1$ and each edge (also called bond) $e \in E$ is independently declared \textit{open} with probability $p$ and \textit{closed} with probability $1-p$. Here $\|x\|_1$ denotes the $\ell_1$-norm of $x \in \mathbb{Z}^d$. The resulting probability measure is denoted by $\mathbb{P}_{p}$.





Let $\left\{ x \leftrightarrow y \right\}$ denote the event that vertices $x$ and $y$ are connected by a finite path of open edges. Let $\mathscr{C}(x)= \left\{ y \in \mathbb{Z}^d \mid x \leftrightarrow y \right\}$ denote the \textit{open cluster} of $x$. It is well known that for $d\geq 2$ there exists a critical probability $p_c \in(0,1)$ for which the model undergoes a phase transition:
\begin{equation}\label{faseovergang}
\mathbb{P}_{p_c}( \exists x \in \mathbb{Z}^d \textit{ s.t. } |\mathscr{C}(x)| = \infty) =  \begin{cases} 0 &  \text{ if } p  < p_c ; \\  1 & \text{  if } p> p_c. \end{cases} 
\end{equation}

Later we will zoom in on what happens at $p=p_c$. Let $\|x\|$ denote the Euclidean norm of $x \in \mathbb{Z}^d$. This choice of norm is not essential, since all norms on $\mathbb{Z}^d$ are equivalent and we only work with estimates that hold up to a constant value. For functions $f$ and $g$, we let $f \asymp g$ denote that $cg \leq f \leq Cg$ holds asymptotically for some constants $c,C>0$. Throughout this article we assume that our lattice is \textit{high-dimensional}, by which we mean that $d>6$  is such that  
\begin{equation}\label{twopointest}\mathbb{P}_{p_c}(x \leftrightarrow y) \asymp \|x- y\|^{2-d},\end{equation}
 for $x,y \in \mathbb{Z}^d$.  
It is widely believed that (\ref{twopointest}) holds in all dimensions $d>6$. In case of nearest-neighbour percolation it has been known for some time that (\ref{twopointest}) is true for all $d\geq 19$ \cite{twopointfunction2} and recently V.d. Hofstad and Fitzner proved it for $d \geq 15$ (in preparation). If there exists an $L>0$ such that $(x,y) \in E \Leftrightarrow \|x-y\| \leq L$, then we speak of \textit{spread-out finite-range} percolation, rather than nearest-neighbour percolation. For this model, it has been proven that (\ref{twopointest}) holds in $d>6$ if the lattice is sufficiently spread out, which means that $L$ should be large enough \cite{twopointfunction}. For readability we restrict ourselves to the case of nearest-neighbour percolation, but all results in this article also hold for spread-out finite-range percolation.

 
In the regime of high dimensions, calculations are relatively easy. In technical practice this is often a consequence of validity of the bound (\ref{twopointest}) on the two-point function, but the intuitive idea behind all this is that for $d$ larger than a certain critical dimension $d_c$, of which the value is believed to be $6$, the model attains \textit{mean-field} behaviour. The amount of space in which open paths can travel has become so large that different pieces of a critical cluster hardly interact. In particular, large open cycles have very small probability. Therefore an open cluster will for many questions behave like a connected graph without cycles: a tree. Percolation on a tree is relatively easy.

\paragraph{Incipient Infinite Cluster}
We now focus on what happens during the phase transition at $p=p_c$. In particular, we want to know how critical clusters behave `as they are becoming infinitely large'. This interpretation is the source of the name \textit{Incipient Infinite Cluster} (IIC), a term originating from the physics literature, which was first defined and treated in a mathematically rigorous way by Kesten \cite{Kesten}. See below for a formal definition.

It turns out that $\mathbb{P}_{p_c}\left(|\mathscr{C}(0)|=\infty \right)=0$ in high dimensions \cite{HaraSladePercolationProb0}, so working directly with $\mathbb{P}_{p_c}$ will not provide us with interesting detailed information about an infinite cluster. This problem can be overcome by conditioning on some event that implies that $|\mathscr{C}(0)|=\infty$, thus  constructing a new probability measure. There exist several constructions of such an IIC-measure that have been proven to be equivalent, providing evidence that the IIC is quite a canonical, robust and unique object. For a precise characterization, the reader is referred to \cite{IICrevisited} and \cite{equivalencePandQmeasure}. We will only directly need the following construction:

\begin{equation}\label{IICconstructie} \mathbb{P}_{\text{IIC}}(F)= \lim_{\|x\| \rightarrow \infty} \mathbb{P}_{p_c}(F \mid  0 \leftrightarrow x) \end{equation}
for cylinder events $F$. In high dimensions, the limit exists irrespective of the direction. 
Through references to literature we will also implicitly use the construction
\begin{equation*}\label{IICconstructie2} \mathbb{Q}_{\text{IIC}}(F)= \lim_{p \uparrow p_c}  \frac{\sum_{x \in \mathbb{Z}^d} \mathbb{P}_p \left( F \cap \left\{ 0 \leftrightarrow x \right\} \right)}{\sum_{x \in \mathbb{Z}^d} \mathbb{P}_{p}\left( 0 \leftrightarrow x \right)}.    \end{equation*}
In high dimensions, the limits $\mathbb{P}_{\text{IIC}}(F)$ and $\mathbb{Q}_{\text{IIC}}(F)$ exist and are equal for all cylinder events $F$. Consequently $\mathbb{P}_{\text{IIC}}$ and $\mathbb{Q}_{\text{IIC}}$ extend to the same probability measure in our context \cite{IICrevisited},\cite{equivalencePandQmeasure}. Expectation value with respect to $\mathbb{P}_{\text{IIC}}$ is denoted by $\mathbb{E}_{\text{IIC}}$. It holds that $\mathbb{P}_{\text{IIC}} \left( |\mathscr{C}(0)|= \infty \right) = 1$ and partly because of this, some authors refer to the IIC as the distribution of $\mathscr{C}(0)$ under $\mathbb{P}_{\text{IIC}}$. However, in the context of $\mathbb{P}_{\text{IIC}}$ the term IIC is also often used to refer to the infinite cluster at the origin itself. We adopt the latter convention.


\begin{Deff*}
IIC is the random graph with vertex set $\mathscr{C}(0)$ and induced edge set
$$\left\{ (x,y) \in \mathscr{C}(0) \times \mathscr{C}(0) \mid (x,y) \text{ is open } \right\}.$$
In many cases we are only interested in the vertices and therefore we abuse notation by writing $\text{IIC}= \mathscr{C}(0)$.
\end{Deff*}

\subsection{Mass dimension and volume growth exponent}

In order to determine how large the (infinite) IIC is, we need to associate some natural notion of dimensionality. On the one hand, we will calculate the mass dimension, which counts the vertices of IIC that are in a cube of finite radius $r$ around the origin. On the other hand, we consider the volume growth exponent, which counts the number of vertices in IIC that can be reached from the origin by an open path of length at most some fixed $r$. In the former case, IIC is counted with respect to the `extrinsic' (Euclidean) metric of the underlying lattice $\mathbb{Z}^d$, while in the latter case, IIC is counted with respect to the `intrinsic' graph distance of the random graph.

\begin{Defs*}

Denote by
$$Q_r = \left\{ x \in \mathbb{Z}^d  \mid \|x\| \leq r \right\}$$
the cube with radius $r$ and boundary
$$\partial Q_r = Q_r \backslash Q_{r-1}.$$
In practice we will want to bound the cardinality of the following three random sets,
$$X_r = \left\{ x \in Q_r \mid 0 \leftrightarrow x \right\}$$
$$X_{r,r} = \left\{ x \in Q_r \mid 0 \overset{Q_r}{\longleftrightarrow} x \right\}$$
$$B_{r}= \left\{ x \in \mathbb{Z}^d   \mid 0 \overset{\leq r}{\longleftrightarrow} x \right\}$$
where $0 \overset{Q_r}{\longleftrightarrow} x$ means that $0$ is connected to $x$ by an open path that does not leave $Q_r$ and $0 \overset{\leq r}{\longleftrightarrow} x$ means that $0$ is connected to $x$ by an open path of length $\leq r$ (with respect to graph distance in the random percolated graph).
\end{Defs*}

\begin{Def*}
The \textit{mass dimension} of a subset $A \subset \mathbb{Z}^d$ is 
$$d_m(A)= \lim_{r \rightarrow \infty} \log_r|A \cap Q_r|$$
if the limit exists.
 The \textit{volume growth exponent} of an infinite connected graph $G$ is defined by
$$d_f(G)= \lim_{r \rightarrow \infty} \log_r |B_G(x,r)|$$
if the limit exists. Here $B_G(x,r)$ is the ball with some center vertex $x$ and radius $r$, with respect to graph distance. 
\end{Def*}
Note that the mass dimension of IIC equals $d_m(\text{IIC})=\lim_{r \rightarrow \infty} \log_{r}|X_r|$ and the volume growth exponent of IIC can be rewritten as $d_f(\text{IIC})=\lim_{r \rightarrow \infty}\log_{r}|B_{\text{IIC}}(0,r)|= \lim_{r \rightarrow \infty} \log_{r}|B_r|$. 

Our main goal is to prove Theorem \ref{MassadimensieIs4}, which  states that on a high-dimensional lattice the mass dimension of IIC almost surely equals $4$.


\begin{Theorem}\label{MassadimensieIs4}
In high dimensions,
$$\mathbb{P}_{\text{IIC}} \left( d_m(IIC) \equiv \lim_{r \rightarrow \infty} \left(\log_{r}|X_r| \right)  = 4 \right) =1.$$
\end{Theorem}

This can be contrasted against Theorem \ref{VolumeGrowthExponentIs2}, which states that on a high-dimensional lattice the volume growth exponent of $\text{IIC}$ almost surely equals $2$. This second result was already implicit in two auxiliary lemmas in \cite{KozmaNachmias}, which we use to obtain a formal derivation of the almost sure statement.
\begin{Theorem}\label{VolumeGrowthExponentIs2}
In high dimensions,
$$\mathbb{P}_{\text{IIC}} \left(d_f(\text{IIC})\equiv \lim_{r \rightarrow \infty} \left(\log_r|B_r|\right) = 2 \right)=1.$$
\end{Theorem}

\subsection{Embedding and conjectures}

\paragraph{On the $4$-dimensionality of IIC.}
Earlier developments in the direction of determining `the' dimension of the IIC include the following. In \cite{equivalencePandQmeasure} it was shown that in high dimensions, $\mathbb{P}_{\text{IIC}}\left( 0 \leftrightarrow x \right) \asymp \|x\|^{4-d}$, implying that $\mathbb{E}_{\text{IIC}}(|X_r|) \asymp C \cdot r^4$. This moment bound, which is also derived in a more general setting in \cite{IICrevisited}, already gave some weak notion of the 4-dimensionality of the IIC. As we will see later, it provides enough information to derive an almost sure upper bound $4$ on the (upper) mass dimension of IIC, essentially using Markov's inequaliy and Borel-Cantelli. However, deriving the corresponding lower bound $4$ on the (lower) mass dimension requires a completely different technique. Concentration inequalities like the second moment method are not powerfull enough \cite{Masterthesis} and many standard techniques from percolation theory don't apply because of the delicate dependency on the origin, induced by the IIC-measure. Indeed, the derivation of the lower bound constitutes the main contribution of this article.


\paragraph{Spectral dimension and other bounds on $|X_r|$ and $|B_r|$.}
The spectral dimension of an infinite connected graph $G$ is defined by
$$d_{s}(G)= -2 \cdot \lim_{r \rightarrow \infty} \log_{r} \left(p_{2r}(x,x) \right)$$
if the limit exists. Here $p_{2r}(x,x)$ is the return probability of a simple random walk on $G$ after $r$ steps. Kozma and Nachmias \cite{KozmaNachmias} showed that $d_s(\text{IIC})= \frac{4}{3}$, thereby confirming the so-called Alexander-Orbach conjecture in high dimensions. For many `nice' graphs and in particular for any Cayley graph $G$ it holds that $d_f(G)=d_s(G)$, but this is not the case for the IIC, as $d_f(\text{IIC})=2 \neq \frac{4}{3}=d_s(\text{IIC})$, suggesting that the IIC is an intrinsically fractal object. Kozma and Nachmias also showed that $\mathbb{E}_{p_c}\left(|B_r| \right) \asymp r$ and $\mathbb{P}_{p_c} \left(B_r \backslash B_{r-1} \neq \emptyset  \right) \asymp r^{-1}.$ These statements are in terms of the intrinsic graph distance and should be contrasted against their extrinsic counterparts $\mathbb{E}_{p_c}\left( |X_r| \right) \asymp r^2$  and $ \mathbb{P}_{p_c} \left(0 \leftrightarrow \partial Q_r \right) \asymp r^{-2}$ \cite{IICrevisited}\cite{KozmaNachmiasArmexponent}.

\paragraph{Growth behaviour of the boundary of $X_{r,r}$.}
In the proof of Theorem \ref{MassadimensieIs4}, we actually also show that $\mathbb{P}_{\text{IIC}}\left( \lim_{r \rightarrow \infty} \log_{r}(|X_{r,r}|)=4 \right)=1$. That is, $|X_{r,r}|$ and $|X_r|$ don't differ very much; they both grow like $r^4$. Define the boundary $\partial X_{r}:= \left\{ x \in \partial Q_r \mid 0 \leftrightarrow x \right\}$. Since $X_r=\bigsqcup_{k=1}^{r} \partial X_{k}$, it is to be expected that $|\partial X_{r}|$ typically grows like $r^3$. Similarly, if we define the `boundary' $\partial X_{k,r}:= \left\{ x \in \partial Q_k \mid 0 \overset{Q_r}{\longleftrightarrow} x \right\}$ then $X_{r,r}=\bigsqcup_{k=1}^{r} \partial X_{k,r}$, so one would expect that $|\partial X_{k,r}|$ grows like $k^3$. We believe this is indeed the case for $k \ll r$, because for those values $|\partial X_{k,r}| \approx |\partial X_{k}|$. However, if $k \approx r$ the picture (presumably) changes completely. Theorem 1.16 in \cite{Masterthesis} yields that there exists a constant $C>0$ such that for all $\lambda, r >0$,   $\mathbb{P}_{\text{IIC}} \left(\sum_{k=1}^{r}|\partial X_{k,k}| \leq \frac{1}{\lambda} \cdot r^3 \right) \leq C \cdot \frac{1}{\lambda}$. A slight adaptation of that proof yields that $\mathbb{P}_{\text{IIC}} \left( |\partial X_{r,r}| \leq \frac{1}{\lambda} \cdot r^2 \right) \leq C \cdot \frac{1}{\lambda}$ and in fact, we conjecture that the opposite bound $\mathbb{P}_{\text{IIC}} \left( |\partial X_{r,r}| \geq \lambda \cdot r^2 \right) \leq C \cdot \frac{1}{\lambda}$ holds too. In other words, we expect $|\partial X_{r,r}|$ to grow like $r^2$ instead of $r^3$. One motivation for the opposite bound comes from Theorem $2$ in \cite{KozmaNachmiasArmexponent}, which essentially says that $|X_{r,r}|$ is smaller than $r^2$ if $|X_r|$ is smaller than $r^4$. To actually prove the opposite bound, it would suffice to show that $\mathbb{E}_{\text{IIC}}\left( |\partial X_{r,r}| \right) \leq C \cdot r^2$, and for this it would be very useful to have a good upper bound on $\mathbb{P}_{\text{IIC}} ( 0 \overset{Q_r}{\longleftrightarrow} x )$, for $x \in \partial Q_r$. While $\mathbb{P}_{\text{IIC}} ( 0 \longleftrightarrow x ) \asymp \|x\|^{4-d}$ depends only on the norm of $\|x\|$ but not really on the choice of norm, the behaviour of $\mathbb{P}_{\text{IIC}} ( 0 \overset{Q_r}{\longleftrightarrow} x )$ is more complicated. For example, if we define the cube $Q_r$ with respect to the $\ell_{\infty}$-norm, then it is much `harder' for an open path that stays entirely inside $Q_r$ to reach a corner vertex $x_1$ of $Q_r$, than it is to reach the center vertex $x_2$ of a face of $Q_r$, although $\|x_1\|_{\infty}=\|x_2\|_{\infty}$.

\paragraph{The backbone of IIC and scaling limits.}
There is a natural subset of the IIC, called the \textit{backbone} (bb) of the IIC, which consists of all open bonds $e=(e_{-},e_{+})$ such that there exist two disjoint open paths, one path from $0$ to $e_{-}$ and the other path from $e_{+}$ to $\infty$. It is expected that the mass dimension of the backbone $\mathbb{P}_{\text{IIC}}$-almost surely equals $2$.
The validity of the almost sure upper bound $2$ is immediate from the known expectation bound $\mathbb{E}_{\text{IIC}}\left( |bb \cap Q_r|\right) \asymp r^2$ \cite{IICrevisited} and an application of Lemma \ref{uppermassbovengrens} from the present article. Heydenreich, V.d. Hofstad, Hulshof and Miermont prepare a proof that the scaling limit of the backbone is a brownian motion, which almost surely has Hausdorff dimension $2$. A related, but wide open conjecture is that the scaling limit of the high-dimensional IIC itself is \textit{Integrated super-Brownian excursion} \cite{scalinglimitIIC}.

\paragraph{The IIC in low dimensions.}
For $d=1$, IIC trivially has mass dimension and volume growth exponent $1$. Kesten proved the bound 
$$ \mathbb{E}_{\text{IIC}} |\text{IIC} \cap Q_r| \asymp r^2 \cdot  \mathbb{P}_{p_c} \left(0 \leftrightarrow \partial Q_r \right),$$ 
which holds for a wide range of lattices on $\mathbb{Z}^2$ \cite{Kesten}. For site percolation on the triangular lattice, Lawler, Schramm and Werner were able to show that $\mathbb{P}_{p_c} \left(0 \leftrightarrow \partial Q_r \right)= r^{-5/48 + o(1)}$  \cite{LSW}. So for this particular lattice, $  \mathbb{E}_{\text{IIC}} |\text{IIC} \cap Q_r|  \asymp r^{2} \cdot r^{-5 /48} = r^{91/48}.$ By the conjectured universality of the exponent, this result presumably holds for all common two-dimensional lattices. Note that $\frac{91}{48}$ is just slightly smaller than $2$, the dimension of the surrounding space.
For $3 \leq d \leq 6$ very little is known rigorously. Simulations by Kumagai suggest that $d_s(\text{IIC})$ ranges from $\approx 1.318 +/-0.001$ for $d=2$ to $\approx 1.34 +/- 0.02$ for $d=5$, which is close to the value $4/3$ that holds in high dimensions, but nevertheless supports the belief that the Alexander-Orbach conjecture is false for $d\leq 6$ \cite{KozmaNachmias}.

\subsection{About the proof}


For Theorem 1 we use an upper bound on the expectation value of $|X_r|$ to derive that $d_m(\text{IIC})\leq 4$, almost surely. The lower bound is the hard (or at least more unusual) part. For this we use the one-arm exponent bound $\mathbb{P}_{p_c}\left(0 \leftrightarrow  \partial Q_r \right) \leq C \cdot \frac{1}{r^2}$, from which it will follow that under $\mathbb{P}_{\text{IIC}}$ a typical shortest open path between $0$ and $\partial Q_r$ has length $r^2$. In Theorem \ref{CrucialeStelling} this is combined with the fact that the intrinsic ball $B_r$ contains approximately $r^2$ elements, yielding that
$|X_r|\geq|X_{r,r}| \approx |B_{(\text{length shortest open path } 0 \leftrightarrow \partial Q_r})| \approx |B_{r^2}| \approx (r^2)^2 = r^4, $ or rather that large downwards deviations of these approximations have small enough probability. The workhorse of this article is Lemma \ref{uppermassbovengrens}, which turns probabilistic bounds into almost sure statements. Indeed, Theorem \ref{VolumeGrowthExponentIs2} follows by a direct application of this lemma to a result from literature.

\section{Ingredients from literature}
In this section we collect ingredients from the literature that we use in our proofs.
\begin{Theorem}[Theorem 1.5 in \cite{IICrevisited}]\label{samenvattingexpectationbounds}
In high dimensions, there exists a constant $C>0$ such that for all $r \geq 1$:
\begin{equation*}\label{claim 4}  \mathbb{E}_{\text{IIC}} \left( |X_r| \right) \leq C \cdot r^4. \end{equation*}
\end{Theorem}


\begin{Theorem}[Corollary of Theorem 1 in \cite{KozmaNachmiasArmexponent}] \label{armexponentBLOE}
In high dimensions, there exists a $C>0$ such that for all $r\geq 1$:
$$\mathbb{P}_{p_c} \left( 0 \leftrightarrow \partial Q_r \right) \leq C \cdot \frac{1}{r^2}.$$
\end{Theorem}

\begin{Lemma}[Lemma 2.5 in \cite{KozmaNachmias}]\label{ReductieTotPc}
In high dimensions, there exists a constant $C>0$ such that for all $r\geq 1$ and any event $E$ measurable with respect to $B_r$ and for any $x \in \mathbb{Z}^d$ with $\|x\|$ sufficiently large:
$$\mathbb{P}_{p_c}(E \cap \left\{0 \leftrightarrow x  \right\}) \leq C \cdot \sqrt{r \cdot \mathbb{P}_{p_c}(E)} \cdot \mathbb{P}_{p_c}(0 \leftrightarrow x).$$
\end{Lemma}

\begin{Lemma}[Essentially Lemma 6.1 in \cite{CycleStructureOfPercolationOnHighDimensionalTori}]\label{kansoprandverbindingmetkleinopenpas}
In high dimensions, there exists a $C>0$ such that for all $\epsilon>0, r\geq 1$:
$$\mathbb{P}_{\text{IIC}} \left(  0  \overset{\leq \epsilon \cdot r^2 }{\longleftrightarrow} \partial Q_r   \right) \leq C \cdot \sqrt{\epsilon},$$
where $\left\{  0  \overset{\leq \epsilon \cdot r^2 }{\longleftrightarrow} \partial Q_r \right\}$ is the event that $0$ is connected to $\partial Q_r$ by an open path of length $\leq \epsilon \cdot r^2$.
\begin{proof}
The event $E= \left\{  0  \overset{\leq \epsilon \cdot r^2 }{\longleftrightarrow} \partial Q_r \right\}$ is measurable with respect to $B_{\epsilon \cdot r^2}$. Therefore, Lemma \ref{ReductieTotPc} implies that for any $x \in \mathbb{Z}^d$ with $\|x\|$ sufficiently large,
$$\mathbb{P}_{p_c} \left(  0  \overset{\leq \epsilon \cdot r^2 }{\longleftrightarrow} \partial Q_r  \mid  0 \leftrightarrow x \right) \leq C' \cdot \sqrt{\epsilon \cdot r^2 \cdot \mathbb{P}_{p_c}\left( 0 \leftrightarrow \partial Q_r \right)} \leq C \cdot \sqrt{\epsilon},$$
where the second inequality follows from Theorem \ref{armexponentBLOE}. Now apply construction (\ref{IICconstructie}) of $\mathbb{P}_{\text{IIC}}$.
\end{proof}
\end{Lemma}

\begin{Lemma}[Essentially Lemmas 2.2 and 2.3 in \cite{KozmaNachmias}]\label{IntrinsicClusterUpperBounds}
In high dimensions, there exists a $C>0$ such that for all $\lambda>1$ and $r \geq 1$:
\begin{equation}\label{idioot1} \mathbb{P}_{\text{IIC}} \left(  |B_{r}| \leq \frac{1}{\lambda} \cdot r^2  \right) \leq C \cdot \frac{1}{\lambda} \end{equation}
and
\begin{equation}\label{idioot2}\mathbb{P}_{\text{IIC}} \left(  |B_{r}| \geq \lambda \cdot r^2  \right) \leq C \cdot \frac{1}{\lambda}. \end{equation}
\begin{proof}
Inequality (\ref{idioot1}) is the statement of Lemma 2.3 in \cite{KozmaNachmias}. On the other hand, Lemma 2.2 in \cite{KozmaNachmias} states that there exists a $C > 0$ such that for all $r \geq 1$ and all $x \in \mathbb{Z}^d$ with $\|x\|$ sufficiently large,
$$\mathbb{E}_{p_c}\left( |B_r| \cdot \mathbbm{1}_{\left\{0 \leftrightarrow x\right\}} \right) \leq C \cdot r^2 \cdot \mathbb{P}_{p_c}\left(0 \leftrightarrow x \right).$$
By Markov's inequality this implies that for all $\lambda>1$ and $r \geq 1$ it holds that $\mathbb{P}_{p_c}\left( |B_r|\geq \lambda\cdot r^2 \mid 0 \leftrightarrow x  \right) \leq C \cdot \frac{1}{\lambda}$, for all $x \in \mathbb{Z}^d$ with $\|x\|$ sufficiently large. Letting $\|x\| \rightarrow \infty$ yields (\ref{idioot2}), because $\left\{ |B_r| \geq \lambda \cdot r^2 \right\}$ is a cylinder event.
\end{proof}
\end{Lemma}

\section{Deriving the main theorems}

The following theorem is crucial for the derivation of Theorem \ref{MassadimensieIs4}. It relies on Lemmas \ref{kansoprandverbindingmetkleinopenpas} and \ref{IntrinsicClusterUpperBounds} and in that sense, it uses that both the cardinality of the intrinsic ball with radius $r$ and the length of the shortest path from $0$ to the boundary of $\partial Q_r$ grow like $r^2$.
\begin{Theorem}\label{CrucialeStelling} 
In high dimensions, there exists a $C>0$ such that for all $\lambda>1$ and $r\geq 1$:
$$\mathbb{P}_{\text{IIC}}\left( |X_{r,r}| \leq \frac{1}{\lambda} \cdot r^4 \right) \leq C \cdot \frac{1}{\lambda^{1/5}}. $$
\begin{proof}
Let $\lambda>1$. Write $\epsilon:= \epsilon(\lambda)= \lambda^{-2/5}$, then
\begin{equation}\label{Splitsing} \mathbb{P}_{\text{IIC}}\left( |X_{r,r}| \leq \frac{1}{\lambda} \cdot r^4 \right) = \mathbb{P}_{\text{IIC}}\left( |X_{r,r}| \leq \frac{1}{\lambda} \cdot r^4 ,  0  \overset{\leq \epsilon \cdot r^2 }{\longleftrightarrow} \partial Q_r \right)  + \mathbb{P}_{\text{IIC}}\left( |X_{r,r}| \leq \frac{1}{\lambda} \cdot r^4 ,  \text{ not } 0  \overset{\leq \epsilon \cdot r^2 }{\longleftrightarrow} \partial Q_r \right). \end{equation}
 By Lemma \ref{kansoprandverbindingmetkleinopenpas} we can bound the first term as follows:
\begin{equation}\label{FirstTerm} \mathbb{P}_{\text{IIC}}\left( |X_{r,r}| \leq \frac{1}{\lambda} \cdot r^4 ,  0  \overset{\leq \epsilon \cdot r^2 }{\longleftrightarrow} \partial Q_r \right) \leq \mathbb{P}_{\text{IIC}}\left( 0  \overset{\leq \epsilon \cdot r^2 }{\longleftrightarrow} \partial Q_r \right) \leq C \cdot \epsilon^{1/2}=C\cdot \frac{1}{\lambda^{1/5}}. \end{equation}
On the other hand, if $\left\{  \text{not } 0  \overset{\leq \epsilon \cdot r^2 }{\longleftrightarrow} \partial Q_r \right\}$ occurs then the intrinsic ball $B_{\epsilon \cdot r^2}$ is a subset of $X_{r,r}$, so $|B_{\epsilon \cdot r^2}| \leq |X_{r,r}|$. Therefore a bound on the second term is given by
\begin{eqnarray}\label{SecondTerm}
\mathbb{P}_{\text{IIC}}\left( |X_{r,r}| \leq \frac{1}{\lambda} \cdot r^4 ,  \text{ not } 0  \overset{\leq \epsilon \cdot r^2 }{\longleftrightarrow} \partial Q_r \right)    &\leq&    \mathbb{P}_{\text{IIC}}\left(  |B_{\epsilon \cdot r^2}| \leq \frac{1}{\lambda} \cdot r^4 \right) \notag \\ 
&=& \mathbb{P}_{\text{IIC}}\left(  |B_{\epsilon \cdot r^2}| \leq \frac{1}{\lambda \cdot\epsilon^2} \cdot (\epsilon \cdot r^2)^2 \right)  \notag \\
&\leq& C \cdot \frac{1}{\lambda \cdot \epsilon^2}  \notag \\
&=& C \cdot \frac{1}{\lambda^{1/5}},
\end{eqnarray}
where the second inequality follows from Lemma \ref{IntrinsicClusterUpperBounds}. Now evaluate (\ref{FirstTerm}) and (\ref{SecondTerm}) in (\ref{Splitsing}) to finish the proof.
\end{proof}
\end{Theorem}

The next lemma will be used to transform the results obtained so far into the almost sure statements of Theorem \ref{MassadimensieIs4} and \ref{VolumeGrowthExponentIs2}. We present a more general and stronger version than	 we actually need.

\begin{Lemma}\label{uppermassbovengrens}
Let $Z_1,Z_2,\ldots$ be a sequence of random variables with values in $\mathbb{R}_{> 0}$, such that $Z_1 \leq Z_2 \leq \ldots$
\begin{enumerate}
\item 
If there exist constants $\beta, \mu, C >0$ such that at least one of the following two conditions holds \begin{itemize} 
\item $\mathbb{E}(Z_r) \leq C \cdot r^{\beta}$ for all $r>0$;
\item $\mathbb{P}(Z_r \geq \lambda \cdot r^{\beta}) \leq C \cdot \frac{1}{\log(\lambda)^{1+\mu}}$ for all $\lambda>1$ and $r >0 $,
\end{itemize}
then:
\begin{equation}\label{groeneharing1}\mathbb{P} \left( \limsup_{r \rightarrow \infty} \left(    \log_{r}(Z_{r}) \right) \leq \beta \right) = 1. \end{equation}
\item
If there exist constants $\alpha, \mu, C >0$ such that at least one of the following two conditions holds \begin{itemize} 
\item $\mathbb{E}\left(\frac{1}{Z_r}\right) \leq C \cdot r^{-\alpha} $ for all $r>0$;
\item $\mathbb{P}(Z_r \leq \frac{1}{\lambda} \cdot r^{\alpha}) \leq C \cdot \frac{1}{\log(\lambda)^{1+\mu}}$ for all $\lambda>1$ and $r > 0 $,
\end{itemize}
then:
\begin{equation}\label{groeneharing2}\mathbb{P} \left( \liminf_{r \rightarrow \infty} \left(    \log_{r}(Z_{r}) \right) \geq \alpha \right) = 1.\end{equation}
\end{enumerate}

\begin{proof}
First note that the first condition of (\ref{groeneharing1})  implies the second condition of (\ref{groeneharing1}). Indeed, by Markov's inequality there exist $C, \mu >0$ such that for all $\lambda >1$ and $r>0$
$$\mathbb{P}(Z_r \geq \lambda \cdot r^{\beta}) \leq \frac{\mathbb{E}(Z_r)}{\lambda \cdot r^{\beta}} \leq \frac{C \cdot r^{\beta}}{\lambda \cdot r^{\beta}} \leq  C \cdot  \frac{1}{\log(\lambda)^{1+\mu}}.$$
Similarly, the first condition of (\ref{groeneharing2})  implies the second condition of (\ref{groeneharing2}). Indeed,
$$\mathbb{P} \left(Z_r \leq \frac{1}{\lambda} \cdot r^{\alpha}\right) = \mathbb{P}\left(\frac{1}{Z_r} \geq \lambda \cdot r^{-\alpha} \right)  \leq \frac{\mathbb{E} \left(\frac{1}{Z_r}\right)}{\lambda \cdot r^{-\alpha}} \leq \frac{C \cdot r^{-\alpha}}{\lambda \cdot r^{-\alpha}} \leq  C \cdot  \frac{1}{\log(\lambda)^{1+\mu}}.$$
It remains to prove (\ref{groeneharing1}) and (\ref{groeneharing2}) under their second condition.

Define the strictly increasing subsequences $r_k=2^k$ and $\lambda_k=2^{ \left(  k^{ \left(  \frac{1+\mu/2}{1+\mu} \right) }  \right) }$. Also define $\epsilon_k:= \log_{r_k}(\lambda_k)=k^{ \left(  \frac{1+\mu/2}{1+\mu} -1 \right) }$. Note that $r_k, \epsilon_k>0$ and $\lambda_k >1$ for all positive integers $k$, and $\lim_{k \rightarrow \infty} \epsilon_{k}=0$. We first prove (\ref{groeneharing1}). For all positive integers $k$ it holds that

\begin{equation}\label{begincentraleafleiding}\mathbb{P} \left( Z_{r_k} \geq \lambda_k \cdot r_k^{\beta} \right) \leq C \cdot \frac{1}{\log(\lambda_k)^{1+\mu}}.\end{equation}
Using the notation $Y_r:= \log_r (Z_r)$ we obtain that
\begin{eqnarray*}
\sum_{k=1}^{\infty} \mathbb{P} \left( Y_{r_k} \geq \epsilon_k +  \beta \right)     &=&     \sum_{k=1}^{\infty} \mathbb{P} \left( Z_{r_k} \geq \lambda_k \cdot r_k^{\beta} \right) \\
&\leq& C \cdot \sum_{k=1}^{\infty} \frac{1}{\log(\lambda_k)^{1+\mu}} \\
&=&  \frac{C}{\log(2)^{1+\mu}} \cdot \sum_{k=1}^{\infty} \frac{1}{k^{1+ \mu / 2}} \\
&<& \infty.
\end{eqnarray*}
By Borel-Cantelli this implies that
$$\mathbb{P} \left( Y_{r_k} \geq \epsilon_k + \beta \text{ for infinitely many } k \right) =0 $$
and because $\lim_{k \rightarrow \infty}\epsilon_k =0$ it follows that
\begin{equation}\label{subsequenceresultaat1} \mathbb{P}\left( \limsup_{k \rightarrow \infty}(Y_{r_k}) \leq \beta \right)=1. \end{equation}
Now consider any $r >0$ and choose $k \in \mathbb{N}$ such that $2^k \leq r \leq 2^{k+1}$. Then
$$ Y_{r} = \frac{\log(Z_{r})}{\log(r)} \leq  \frac{\log(Z_{2^{k+1}})}{\log(2^k)} =   \frac{\log(Z_{2^{k+1}})}{\log(2^{k+1})} \cdot \frac{\log(2^{k+1})}{\log(2^k)} = Y_{2^{k+1}} \cdot \frac{k+1}{k}$$
and
$$  Y_{r} = \frac{\log(Z_{r})}{\log(r)} \geq  \frac{\log(Z_{2^{k}})}{\log(2^{k+1})} =   \frac{\log(Z_{2^{k}})}{\log(2^{k})} \cdot \frac{\log(2^{k})}{\log(2^{k+1})} = Y_{2^{k}} \cdot \frac{k}{k+1},$$
so 
\begin{equation}\label{limsupYsigmak1} \limsup_{r\rightarrow \infty} Y_{r} = \limsup_{k \rightarrow \infty} Y_{2^{k}}\end{equation}
and
\begin{equation}\label{liminfYsigmak1} \liminf_{r\rightarrow \infty} Y_{r} = \liminf_{k \rightarrow \infty} Y_{2^{k}}.\end{equation}
Evaluating (\ref{limsupYsigmak1}) in (\ref{subsequenceresultaat1}) yields the desired statement (\ref{groeneharing1}).

The proof of (\ref{groeneharing2}) is almost the same. By the arguments used in (\ref{begincentraleafleiding}) - (\ref{subsequenceresultaat1}) we obtain 
$$\mathbb{P} \left( Y_{r_k} \leq - \epsilon_k + \alpha \text{ for infinitely many } k \right) =0 $$
and therefore
\begin{equation}\label{subsequenceresultaat2} \mathbb{P}\left( \liminf_{k \rightarrow \infty}(Y_{r_k}) \geq \alpha \right)=1. \end{equation}

Evaluating (\ref{liminfYsigmak1}) in (\ref{subsequenceresultaat2}) yields the desired statement (\ref{groeneharing2}).\\
\end{proof}
\end{Lemma}
We are ready to prove the main theorems.\\

\textit{Proof of Theorem \ref{MassadimensieIs4}}.\\
Apply Lemma \ref{uppermassbovengrens}.(i) to Theorem \ref{samenvattingexpectationbounds}, with $Z_r=|X_r|$ and $\beta=4$, to obtain
\begin{equation} \label{pensioenakkoord1} \mathbb{P}_{\text{IIC}} \left( \limsup_{r \rightarrow \infty} \left(\log_{r}|X_r| \right) \leq 4 \right)=1. \end{equation}
Apply Lemma \ref{uppermassbovengrens}.(ii) to Theorem \ref{CrucialeStelling}, with $Z_r=|X_{r,r}|$ and $\alpha=4$, to obtain
\begin{equation}\label{pensioenakkoord2}\mathbb{P}_{\text{IIC}} \left( \liminf_{r \rightarrow \infty} \left(\log_{r}|X_{r,r}| \right) \geq 4 \right)=1. \end{equation}
Because $|X_{r,r}| \leq |X_r|$ for all $r \geq 0$ the theorem now follows from (\ref{pensioenakkoord1}) and (\ref{pensioenakkoord2}).\\
\qed

\textit{Proof of Theorem \ref{VolumeGrowthExponentIs2}}.\\
Apply Lemma \ref{uppermassbovengrens}.(i) and \ref{uppermassbovengrens}.(ii) to Lemma \ref{IntrinsicClusterUpperBounds}, with $Z_r=|B_r|$ and $\alpha=\beta=2$, to obtain
$$\mathbb{P}_{\text{IIC}} \left( \limsup_{r \rightarrow \infty} \left(\log_{r}|B_r| \right) \leq 2 \right)=\mathbb{P}_{\text{IIC}} \left( \liminf_{r \rightarrow \infty} \left(\log_{r}|B_r| \right) \geq 2 \right)=1.$$
\qed



\end{document}